\documentclass[12pt]{article}
\usepackage{amssymb,amsthm,amscd,array}

\let\ssection=\section
\renewcommand{\section}{\setcounter{equation}{0}\ssection}
\newcommand{\bbR}{\mathbb{R}}
\newcommand{\bbZ}{\mathbb{Z}}

\newcommand{\bbC}{\mathbb{C}}

\newcommand{\ad}{\mathrm{ad}}

\newcommand{\cD}{{\mathcal{D}}}

\newcommand{\Div}{\mathrm{Div}}

\newcommand{\End}{\mathrm{End}}
\newcommand{\cF}{{\mathcal{F}}}

\newcommand{\GL}{{\mathrm{GL}}}

\newcommand{\gr}{\mathrm{gr}}

\newcommand{\Hom}{\mathrm{Hom}}

\newcommand{\cL}{{\mathcal{L}}}

\newcommand{\cR}{{\mathcal{R}}}

\newcommand{\cS}{{\mathcal{S}}}

\newcommand{\Vect}{\mathrm{Vect}}

\newcommand{\half}{\frac{1}{2}}
\newcommand{\fg}{\mathfrak{g}}

\begin{document}

%\baselineskip=18pt

%\textwidth=16truecm
%\textheight=24truecm
%\hoffset=-1.5truecm
%\voffset=-2.5truecm

\def\a{\alpha}
\def\b{\beta}
\def\d{\delta}
\def\g{\gamma}
\def\om{\omega}
\def\r{\rho}
\def\s{\sigma}
\def\vf{\varphi}
\def\l{\lambda}
\def\m{\mu}
\def\k{\kappa}

\def\implies{\Rightarrow}

\oddsidemargin .1truein
\newtheorem{thm}{Theorem}[section]
\newtheorem{lem}[thm]{Lemma}
\newtheorem{cor}[thm]{Corollary}
\newtheorem{pro}[thm]{Proposition}
\newtheorem{ex}[thm]{Example}
\newtheorem{rmk}[thm]{Remark}
\newtheorem{defi}[thm]{Definition}
%\newremark{ex}[thm]{Example}

%\newtheorem{thm}{Theorem}
%\newtheorem{lem}{Lemma}
%\newtheorem{cor}{Corollary}
%\newtheorem{prop}[thm]{Proposition}
%\newtheorem{definition}{Definition}

\title{Multi-parameter deformations of the
module of symbols of differential operators}

\author{
B.~Agrebaoui\footnote{Facult\'e des Sciences de Sfax, BP. 802 3018 Sfax Tunisie,
B.Agreba@fss.rnu.tn}\and
F.~Ammar\footnote{Facult\'e des Sciences de Sfax, BP. 802 3018 Sfax Tunisie,
Faouzi.Ammar@fss.rnu.tn}
\and
V.~Ovsienko\footnote{CNRS, Luminy Case 907,
F--13288 Marseille, Cedex 9, France
ovsienko@cpt.univ-mrs.fr}
}

\date{}

\maketitle

\thispagestyle{empty}

\begin{abstract}
The space of symbols of differential operators on a smooth manifold
(i.e., the space of symmetric contravariant tensor fields)
is naturally a module over the Lie algebra of vector fields. 
We study, in the case of $\bbR^n$ with $n\geq2$, 
multi-parameter formal deformations of this module. 
The space of linear differential operators
on $\bbR^n$ provides an important class
of such formal deformations; we show, however, 
that the whole space of deformations is much larger.
\end{abstract}

\vskip1cm
\noindent
%\textbf{Keywords:} Deformations, modules of differential
%operators, cohomology of Lie algebras of vector fields.

\newpage

%%%%%%%%%%%%%%%%%%%%%%%%%%%%%%%%%%%%%%%%%%%%%%%%%%%%%%%%%%%%%%%%%%%%%%%%%%%%%%%%%%%%
%%%%%%%%%%%%%%%%%%%%%%%%%%%%%%%%%%%%%%%%%%%%%%%%%%%%%%%%%%%%%%%%%%%%%%%%%%%%%%%%%%%%
\section{Introduction}
%%%%%%%%%%%%%%%%%%%%%%%%%%%%%%%%%%%%%%%%%%%%%%%%%%%%%%%%%%%%%%%%%%%%%%%%%%%%%%%%%%%%
%%%%%%%%%%%%%%%%%%%%%%%%%%%%%%%%%%%%%%%%%%%%%%%%%%%%%%%%%%%%%%%%%%%%%%%%%%%%%%%%%%%%

The space of linear differential operators on tensor densities
over a smooth manifold is naturally a module over the Lie algebra of vector fields. 
This module structure has been studied in a series of recent papers 
(see \cite{DO,LO1,CMZ,LO2,Gar,MAT,BO,CM} and references therein).
The module of differential operators can be viewed
as a deformation of the corresponding module of symbols;
the general framework of the deformation theory (see e.g. \cite{G,NR,R,F,FF}),
therefore, relates its study to the
cohomology of the Lie algebra of vector fields, cf. \cite{DO,LO2}.

The first cohomology space of the Lie algebra of vector fields, classifying the
infinitesimal deformations of the module of symbols has been calculated,
for an arbitrary smooth manifold, in \cite{LO2}
(see also \cite{BO} for the details in the one-dimensional case).
Of course, not for every infinitesimal deformation there exists a formal
deformation containing the latter as an infinitesimal part.
The obstructions are characterized
in terms of Nijenhuis-Richardson products of non-trivial
first cohomology classes. The main problem considered in this paper is
to determine the integrability condition,
i.e., a necessary and sufficient condition for an infinitesimal deformation
that guarantees existence of a formal deformation.
We provide such a condition in the case of $\bbR^n$, where $n\geq2$.

Let $\cF_\l$ be the space of tensor densities of degree $\l\in\bbR$ on $\bbR^n$.
The two-parameter family of $\Vect(\bbR^n)$-modules, $\cD_{\l,\m}$
of linear differential operators from $\cF_\l$ to $\cF_\m$
will provide us with an important class of examples of non-trivial deformations
of the module of symbols.

The classical deformation theory traditionally deals with one-parameter 
deformations (cf. \cite{G,NR,R}). We will study multi-prarameter deformations
and adopt here a modern viewpoint of miniversal deformations 
(see \cite{FF}). Our methods are similar to those of
\cite{OR,OR1}.
The first cohomology space of $\Vect(\bbR^n)$ has in our case
a canonical basis; we consider a commutative algebra
generated by the parameters of deformation, corresponding to 
all non-trivial cohomology classes. This allows us to consider the
most general multi-prarameter infinitesimal deformation.
The obstructions for integrability of an infinitesimal deformation 
is expressed in terms of algebraic relations between the generators.

%\medskip

\textbf{Acknowledgments:} 
We are grateful to P.~Lecomte for fruitfull discussions that
considerably simplified our proofs and for careful reading of a preliminary
version of this paper and also to C.~Duval and C.~Roger for numerous
enlightening discussions. The third author thanks le service de G\'eom\'etrie et
Th\'eorie des Algorithmes de L'Universit\'e de Li\`ege where a part of this work
was done  and F.~Boniver for help.

\goodbreak

%%%%%%%%%%%%%%%%%%%%%%%%%%%%%%%%%%%%%%%%%%%%%%%%%%%%%%%%%%%%%%%%%%%%%%%%%%%%%%%%%%%%
%%%%%%%%%%%%%%%%%%%%%%%%%%%%%%%%%%%%%%%%%%%%%%%%%%%%%%%%%%%%%%%%%%%%%%%%%%%%%%%%%%%%
\section{The general framework}\label{GFW}
%%%%%%%%%%%%%%%%%%%%%%%%%%%%%%%%%%%%%%%%%%%%%%%%%%%%%%%%%%%%%%%%%%%%%%%%%%%%%%%%%%%%
%%%%%%%%%%%%%%%%%%%%%%%%%%%%%%%%%%%%%%%%%%%%%%%%%%%%%%%%%%%%%%%%%%%%%%%%%%%%%%%%%%%%

Let us start with the notion of (multi-parameter) deformations
over a commutative algebra.
Our approach will be similar to those of \cite{OR,OR1};
it corresponds to the notion of miniversal
deformations \cite{FF} in a special case when one can choose a basis
of the first cohomology space.

%%%%%%%%%%%%%%%%%%%%%%%%%%%%%%%%%%%%%%%%%%%%%%%%%%%%%%%%%%%%%%%%%%%%%%%%%%%%%%%%%%%%
\subsection{Polynomial deformations}
%%%%%%%%%%%%%%%%%%%%%%%%%%%%%%%%%%%%%%%%%%%%%%%%%%%%%%%%%%%%%%%%%%%%%%%%%%%%%%%%%%%%

Let $\fg$ be a Lie algebra and $(V,\rho)$ a $\fg$-module, where $V$ is a vector
space and $\rho$ is a homomorphism $\rho:\fg\to\End(V)$. We will consider
multi-paraleter formal deformations, i.e., formal series 
\begin{equation}
\label{def}
\r(t)=\r+ 
\sum_{1\leq m<\infty}\vf_m(t)
\end{equation}
where $t=(t_1,\ldots,t_p)$ are the parameters of deformation
and each term $\vf_m(t)$ is a linear map 
$
\vf_m(t):\fg\to
\End(V)\otimes\bbC[t]
$
which is a
homogeneous polynomial in $t$ of degree~$m$.
The expression $\r(t)$ must satisfy the homomorphism
condition, that is, for every $X,Y\in\fg$
\begin{equation}
\label{hom}
\r(t)([X,Y])
=
[\r(t)(X),\r(t)(Y)]
\end{equation}
where the bracket in the right hand side stands for the standard commutator
in~$\End(V)$ extended to the formal series $\End(V)\otimes\bbC[[t]]$.

%\medskip
%\noindent
%{\bf Remark}.
%One may consider $t_i$ as formal parameters and speak about formal deformations, or
%as real (or complex) numbers. We will use in the sequel the both possibilities.

%%%%%%%%%%%%%%%%%%%%%%%%%%%%%%%%%%%%%%%%%%%%%%%%%%%%%%%%%%%%%%%%%%%%%%%%%%%%%%%%%%%%
\subsection{The Maurer-Cartan equation}
%%%%%%%%%%%%%%%%%%%%%%%%%%%%%%%%%%%%%%%%%%%%%%%%%%%%%%%%%%%%%%%%%%%%%%%%%%%%%%%%%%%%

The standard Chevalley-Eilenberg differential (see \cite{F}) retains,
in the case of linear maps from $\fg$ to $\End(V)$ to the following formula.
Given a linear map $a:\fg\to\End(V)$, its differential $\d{}a$ is the bilinear 
skew-symmetric map 
$$
\d{}a(X,Y)
=
a([X,Y])-[\r(X),a(Y)]+[\r(Y),a(X)]
$$
The standard cup-product of linear maps
$a,b:\fg\to\End(V)$ is a bilinear map 
$[\![a,b]\!]:\fg\otimes\fg\to\End(V)$ defined by
\begin{equation}
\label{cup}
[\![a,b]\!](X,Y)
=
-[a(X),b(Y)]+[a(Y),b(X)]
\end{equation}
It is also called the Nijenhuis-Richardson product \cite{NR}.

Put $\vf(t)=\r(t)-\r$, one easily checks that the condition (\ref{hom}) reads
\begin{equation}
\label{MC}
\d\vf(t)
+
\half[\![\vf(t),\vf(t)]\!] =0
\end{equation}
This is the Maurer-Cartan equation (also called the deformation equation,
cf. \cite{NR}). Although it is equivalent
to~(\ref{hom}), it is usefull to relate the deformations (\ref{def}) with
the cohomology theory.

%%%%%%%%%%%%%%%%%%%%%%%%%%%%%%%%%%%%%%%%%%%%%%%%%%%%%%%%%%%%%%%%%%%%%%%%%%%%%%%%%%%%
\subsection{Equivalent deformations}
%%%%%%%%%%%%%%%%%%%%%%%%%%%%%%%%%%%%%%%%%%%%%%%%%%%%%%%%%%%%%%%%%%%%%%%%%%%%%%%%%%%%

Two deformations $\r(t)$ and $\r'(t)$ are called {\it equivalent} if there exists
an inner automorphism $I(t):\End(V)\otimes\bbC[[t]]\to\End(V)\otimes\bbC[[t]]$ 
of the form
\begin{equation}
\label{autom}
I(t)=\exp\Big(
\sum_{1\leq{}i\leq{}p}t_i\,\ad{}A_i
+\sum_{1\leq{}i,j\leq{}p}t_it_j\,\ad{}A_{ij}
+\cdots
\Big),
\end{equation}
where $A_i, A_{ij},\ldots$ are some elements of $\End(V)$,
satisfying the relation
\begin{equation}
\label{equiv}
I(t)\r(t)=\r'(t).
\end{equation}

%%%%%%%%%%%%%%%%%%%%%%%%%%%%%%%%%%%%%%%%%%%%%%%%%%%%%%%%%%%%%%%%%%%%%%%%%%%%%%%%%%%%
\subsection{Infinitesimal deformations and first cohomology}
%%%%%%%%%%%%%%%%%%%%%%%%%%%%%%%%%%%%%%%%%%%%%%%%%%%%%%%%%%%%%%%%%%%%%%%%%%%%%%%%%%%%

The first-order term $\vf_1(t)$ of the expression (\ref{def}) is called an
infinitesimal deformation. It is of the form
\begin{equation}
\label{inf}
\vf_1(t)=t_1c_1+\cdots +t_pc_p.
\end{equation}

It is easy to check that the equation (\ref{MC}) implies that each linear map
$c_i:\fg\to\End(V)$ is a 1-cocycle (see e.g. \cite{F}) for the details).  
Furthermore, if $\r(t)$ and $\r'(t)$ are equivalent deformations, then
the corresponding cocycles in the infinitesimal
deformations are cohomologous, namely 
$c_i=c_i'+\d{}A_i$. Therefore,
an infinitesimal deformation is defined, up to equivalence, by
the cohomology classes
$[c_1],\ldots,[c_p]$ in $H^1(\fg;\End(V))$.

It is natural to assume the classes
$[c_1],\ldots,[c_p]$ linearly independent. Moreover, we will choose a basis of
$H^1(\fg;\End(V))$ and consider the most general multi-parameter
deformation. 

%%%%%%%%%%%%%%%%%%%%%%%%%%%%%%%%%%%%%%%%%%%%%%%%%%%%%%%%%%%%%%%%%%%%%%%%%%%%%%%%%%%%
\subsection{Obstructions and commutative algebras}\label{Obstr}
%%%%%%%%%%%%%%%%%%%%%%%%%%%%%%%%%%%%%%%%%%%%%%%%%%%%%%%%%%%%%%%%%%%%%%%%%%%%%%%%%%%%

Given an infinitesimal deformation of a $\fg$-module $V$, it is called integrable
if there exists a formal deformation containing it as an infinitesimal part.
%In other words, it can be prolonged to a formal deformation.

Developing (\ref{MC}), one obtains a recurrent system
\begin{equation}
\label{ordk}
\d\vf_m(t)+\half\sum_{i+j=m}[\![\vf_i(t),\vf_j(t)]\!]=0
\end{equation}

The second-order term in (\ref{ordk}) is
$
\d\vf_2(t)+\half[\![\vf_1(t),\vf_1(t)]\!]=0.
$
The cohomology class of $[\![\vf_1(t),\vf_1(t)]\!]$
is, therefore, an obstruction to existence of the second-order term 
$\vf_2(t)$.
It is an element of $H^2(\fg;\End(V))\otimes\bbC[t]$, where
polynomial coefficients that are homogeneous second-order polynomials in $t$.
For existence of $\vf_2(t)$ it is necessary and sufficient that
these obstructions vanish.
One thus obtains second-order relations for the
parameters 
$t_1,\ldots,t_p$. 

In the same way, each term of the system (\ref{ordk})
is a homogeneous polynomial  of order $m$ in $t$.
This leads to a system of algebraic relations
on the formal parameters: $R_m(t)=0$, where $m\geq2$. 
To construct multi-parameter formal deformations of the form (\ref{def}),
one has to consider a commutative associative algebra
generated by $t_1,\cdots,t_p$ such that
all the relations, $R_m(t)=0$ are satisfied.

A notion of versal deformation introduced in \cite{FF} is a universal object
of the category of multi-parameter deformations. Any multi-parameter deformation
can be obtained from the versal deformation as a homomorphism of the corresponding
commutative algebras. If one chooses a basis 
$[c_1],\ldots,[c_p]$ in $H^1(\fg;\End(V))$, the versal deformation corresponds
to the commutative algebra 
$$
A=\bbC[t_1,\ldots,t_p]/R,
$$
where $R$ is the ideal generated by the relations, $R_m(t)=0$.

It worth noticing that up to now we were considering only deformations
with a finite number of parameters (just as in the above definitions).
However, following \cite{FF}, we will include into the considerations the case 
of graded modules with infinitely many independent parameters of deformation.

%%%%%%%%%%%%%%%%%%%%%%%%%%%%%%%%%%%%%%%%%%%%%%%%%%%%%%%%%%%%%%%%%%%%%%%%%%%%%%%%%%%%
%%%%%%%%%%%%%%%%%%%%%%%%%%%%%%%%%%%%%%%%%%%%%%%%%%%%%%%%%%%%%%%%%%%%%%%%%%%%%%%%%%%%
\section{Deformations of $\bbZ$-graded modules}\label{defGrad}
%%%%%%%%%%%%%%%%%%%%%%%%%%%%%%%%%%%%%%%%%%%%%%%%%%%%%%%%%%%%%%%%%%%%%%%%%%%%%%%%%%%%
%%%%%%%%%%%%%%%%%%%%%%%%%%%%%%%%%%%%%%%%%%%%%%%%%%%%%%%%%%%%%%%%%%%%%%%%%%%%%%%%%%%%

Let us consider in more details a particular case when the module
$(V,\r)$ is split into a direct sum of $\fg$-modules:
\begin{equation}
\label{decompMod}
V=\bigoplus_{k\in\bbZ}V_k
\end{equation}
Suppose that for some values $i\in\bbZ$ there exist
non-trivial cocycles $c_i$ on $\fg$ with values in $\End(V)$ such that for all
$X\in\fg$ one has 
\begin{equation}
\label{HomCoc}
c_i(X)|_{V_k}\subset{}V_{k-i}.
\end{equation}
Assume, furthermore, that there is a deformation of the form:
\begin{equation}
\label{oldAct}
\r(\tau)=
\r+\sum_{i\in\bbZ}\tau_i\,c_i+(\tau^2)
\end{equation}
where $\tau_i$ are the \textit{free} parameters, i.e., 
the paremeters generate the free commutative algebra
$\bbC[\tau_i]$.
%Note that one can deal here with a formal deformation
%as well as with a ``genuine'' deformation
%depending on real (or complex) parameters $\tau_i\in\bbR$ (or $\bbC$).

The following construction is meant to use the extra degrees of freedom
given by the decomposition (\ref{decompMod}). We will add formal parameters
indexed by $k\in\bbZ$.
Consider each cocycle
$c_i^k:\fg\to\Hom(V_k,V_{k-i})$ defined by the restriction:
\begin{equation}
\label{restCoc}
c_i^k(X):=\left.c_i(X)\right|_{V_k}
\end{equation}
as independent. 

\begin{pro}
\label{proGen}
There exists a formal deformation of the form
\label{mainDef}
\begin{equation}
\label{decompAct}
\r(t)=
\r+\sum_{i,k\in\bbZ}\,t_i^k\,c_i^k+(t^2)
\end{equation}
where $t_i^k$ are formal parameters satisfying the relation
\begin{equation}
\label{relT}
t_i^{k-j}\,{}t_j^k=t_i^k\,{}t_j^{k-i}.
\end{equation}
\end{pro}
\begin{proof}
The original deformation (\ref{oldAct}) satisfies the Maurer-Cartan equation
(\ref{MC}). In each order $m$ the equation
(\ref{ordk}) for the deformation (\ref{oldAct})
has a solution $\r_m(\tau)\in\Hom(\fg;\End(V))\otimes\bbC[\tau]$ which is a
homogeneous polynomial in $\tau$ of degree $m$.
Replacing in $\r_m(\tau)|_{V_k}$
each monomial $\tau_{i_1}\cdots\tau_{i_{m-1}}\tau_{i_m}$ by
$t_{i_1}^{k-i_2-\cdots-i_m}\cdots{}t_{i_{m-1}}^{k-i_m}t_{i_m}^k$,
one, obviously, gets a solution $\r_m(t)$.
\end{proof}

%\medskip
%\noindent
%{\bf Remark}.
%Proposition \ref{proGen} means that in the case when the parameters
%$\tau_i$ are free, the relations (\ref{relT}) are sufficient for existence of a
%formal deformation (\ref{decompAct}), or, in onther words, for integrability
%of the corresponding infinitesimal deformation. If the parameters
%$\tau_i$ satisfy some relations $R(\tau)$, one readily writes the corresponding
%relations $R^k(t)$ for each $k$ that together with (\ref{relT}) will be sufficient
%integrability conditions.

%%%%%%%%%%%%%%%%%%%%%%%%%%%%%%%%%%%%%%%%%%%%%%%%%%%%%%%%%%%%%%%%%%%%%%%%%%%%%%%%%%%%
%%%%%%%%%%%%%%%%%%%%%%%%%%%%%%%%%%%%%%%%%%%%%%%%%%%%%%%%%%%%%%%%%%%%%%%%%%%%%%%%%%%%
\section{The main results}
%%%%%%%%%%%%%%%%%%%%%%%%%%%%%%%%%%%%%%%%%%%%%%%%%%%%%%%%%%%%%%%%%%%%%%%%%%%%%%%%%%%%
%%%%%%%%%%%%%%%%%%%%%%%%%%%%%%%%%%%%%%%%%%%%%%%%%%%%%%%%%%%%%%%%%%%%%%%%%%%%%%%%%%%%

%%%%%%%%%%%%%%%%%%%%%%%%%%%%%%%%%%%%%%%%%%%%%%%%%%%%%%%%%%%%%%%%%%%%%%%%%%%%%%%%%%%%
\subsection{The space of symbols}
%%%%%%%%%%%%%%%%%%%%%%%%%%%%%%%%%%%%%%%%%%%%%%%%%%%%%%%%%%%%%%%%%%%%%%%%%%%%%%%%%%%%

Consider the Lie algebra $\Vect(\bbR^n)$ of smooth vector fields on $\bbR^n$
and the space $\cS$ of smooth symmetric contravariant tensor
fields on $\bbR^n$. 
The space $\cS$ is naturally isomorphic to the space
of functions on $T^*\bbR^n$ polynomial on fibers. Clearly, $\cS$ has
a structure of a Poisson algebra with natural graduation
\begin{equation}
\label{decomp}
\cS=\bigoplus_{k=0}^\infty\cS_k,
\end{equation}
where $\cS_k$ is the space of $k$-th order tensor fields.

The space $\cS$ is a $\Vect(\bbR^n)$-module since
$\Vect(\bbR^n)\subset\cS$.
In Darboux coordinates,
the action of $X\in\Vect(\bbR^n)$ on $\cS$ is given by the Hamiltonian vector 
field\footnote{here and below sum over repeated indices is understood.}
\begin{equation}
\label{Hamilton}
L_X=
\frac{\partial{}X}{\partial\xi_i}\,
\frac{\partial}{\partial{}x^i}
-
\frac{\partial{}X}{\partial{}x^i}\,
\frac{\partial}{\partial\xi_i}\,,
\end{equation}
which is nothing but the Lie derivative of tensor fields.

The aim of this paper is to study multi-parameter formal
deformations of this module. We will restrict our considerations to the 
multi-parameter formal deformations which are {\it differentiable},
i.e., each term in the formal series (\ref{def}) supposed to be a differential
operator on $\cS$.

%%%%%%%%%%%%%%%%%%%%%%%%%%%%%%%%%%%%%%%%%%%%%%%%%%%%%%%%%%%%%%%%%%%%%%%%%%%%%%%%%%%%
\subsection{Description of the infinitesimal deformations}
\label{CohMulti}
%%%%%%%%%%%%%%%%%%%%%%%%%%%%%%%%%%%%%%%%%%%%%%%%%%%%%%%%%%%%%%%%%%%%%%%%%%%%%%%%%%%%

According to the general framework, one needs an information about the space of
the first cohomology of $\Vect(\bbR^n)$ with coefficients in~$\End(\cS)$
in order to describe the infinitesimal deformations. The module
$\End(\cS)$ is decomposed as follows:
$$
\End(\cS)=\bigoplus_{k,\ell}\Hom(\cS_k,\cS_\ell).
$$
To study the $\Vect(\bbR^n)$-cohomology with coefficients in $\End(\cS)$
it then suffice to consider the cohomology with
coefficients in each module $\Hom(\cS_k,\cS_\ell)$.
We will, furthermore, restrict ourself to the subspace
$\cD(\cS_k,\cS_\ell)\subset\Hom(\cS_k,\cS_\ell)$ given by
differential operators from $\cS_k$ to $\cS_\ell$.

The space of first cohomology of the Lie algebra of vector fields with
coefficients in~$\cD(\cS_k,\cS_\ell))$ has been calculated, for an arbitrary manifold
$M$ of $\dim{M}\geq2$, in~\cite{LO2}. We recall here the result
in the case $M=\bbR^n$.
\begin{equation}
\label{CohomResult}
H^1(\Vect(\bbR^n);\cD(\cS_k,\cS_\ell))
=
\left\{
\begin{array}{lcl}
\bbR & , & \hbox{if} \quad
\left[
\begin{array}{l}
k-\ell=0,\\
k-\ell = 1 \hbox{ and } \ell\neq0,\\
k-\ell = 2 
\end{array}
\right.
\\[16pt]
0 & , & \hbox{otherwise}
\end{array}
\right.
\end{equation}
One has, therefore, infinitely many non-trivial cohomology classes 
generating an infinitesimal deformation of the $\Vect(\bbR^n)$-module $\cS$.

Let us give the explicit formul{\ae} for corresponding 1-cocycles.

\medskip

a) For all $k\geq0$ there is a 1-cocycle with values in $\cD(\cS_k,\cS_k)$ that
associates to $X\in\Vect(\bbR^n)$ the operator of multiplication by the function
\begin{equation}
\label{Cocycle0}
c_0(X)=
\Div(X)
\end{equation}
%although this formula does not depend on $k$,
%we use the notation $c^k_0$ to distinguish the cocycles in each order $k$.

b) For all $k\geq2$ there is a 1-cocycle with values in $\cD(\cS_k,\cS_{k-1})$
given by
\begin{equation}
\label{Cocycle1}
c_1(X)=
\frac{\partial^2 X}{\partial x^i\partial x^j}\,
\frac{\partial^2}{\partial \xi_i\partial \xi_j}
\end{equation}

\medskip

\noindent
{\bf Remark}.
More geometrically, this cocycle can be written as the Lie derivative
of the (flat) connection on $\bbR^n$, namely,
$c_1(X)=L_X(\nabla)$.

\medskip

c) For all $k\geq2$ there is a 1-cocycle with values in $\cD(\cS_k,\cS_{k-2})$
given by
\begin{equation}
\label{Cocycle2}
c_2(X)=
\frac{\partial^3 X}{\partial x^i\partial x^j\partial x^l}\,
\frac{\partial^3}{\partial \xi_i\partial \xi_j\partial \xi_l}
-
3\,\frac{\partial^3 X}{\partial x^i\partial x^j\partial \xi_l}\,
\frac{\partial^2}{\partial \xi_i\partial \xi_j}
\frac{\partial}{\partial x^l}
\end{equation}

\medskip

\noindent
{\bf Remark}.
This cocycle is related to the
famous Moyal product, namely for $P\in\cS_k$,
$c_2(X)(P)$ coincides with the trird-order term in the Moyal product
of $X$ and~$P$.

\medskip

As in Section \ref{defGrad}, we will use the notation
$$
c_i^k={c_i}\left|_{\cS_k}\right.,
\qquad
i=0,1,2
$$
and deal with independent cocycles
$c^k_0,c^k_1,c^k_2$.

%%%%%%%%%%%%%%%%%%%%%%%%%%%%%%%%%%%%%%%%%%%%%%%%%%%%%%%%%%%%%%%%%%%%%%%%%%%%%%%%%%%%
\subsection{Integrability condition}
%%%%%%%%%%%%%%%%%%%%%%%%%%%%%%%%%%%%%%%%%%%%%%%%%%%%%%%%%%%%%%%%%%%%%%%%%%%%%%%%%%%%

According to the results of \cite{LO2} (see Section \ref{CohMulti}),
the infinitesimal deformations of the Lie derivative (\ref{Hamilton})
are of the form: $\r(t)(X)=L_X+\vf_1(t)(X)$ with
\begin{equation}
\label{infMulti}
\vf_1(t)=
\sum_{0\leq{}k<\infty}
t_0^k\,c^k_0 + 
\sum_{2\leq{}k<\infty}
\left(t_1^k\,c^k_1 + t_2^k\,c^k_2\right)
\end{equation}
where the symbols $t_0^k,t_1^k,t_2^k$ stand for independent formal parameters.
We, therefore, have to deal with infinitesimal deformations with infinite
number of parameters.

Let us formulate the main result of this paper.
\begin{thm}
\label{multiThm}
The following relations

\noindent
a) one series of second-order relation $R^k_2(t)$:
\begin{equation}
\label{condMulti1}
t_1^k\,t_2^{k-1} - t_1^{k-2}\,t_2^k =0
\quad,\qquad
k\geq4
\end{equation}

\noindent
b) two series of third-order relations, namely $R^k_3(t)$:
\begin{equation}
\label{condMulti2}
\quad
\left(t_0^k-t_0^{k-1}\right)
t_1^k\,t_2^{k-1}
=0
\quad,\qquad
k\geq3
\end{equation}
and $\widetilde{R}^k_3(t)$:
\begin{equation}
\label{condMulti3}
\quad
\left(t_0^k-t_0^{k-2}\right)
t_2^k\,t_2^{k-2} = 0
\quad,\qquad
k\geq4
\end{equation}
are necessary and sufficient for 
integrability of the infinitesimal deformation (\ref{infMulti}).
\end{thm}

\noindent
The following statement can be considered as a corollary of Theorem
\ref{multiThm}, but we will give its elementary proof.

\begin{pro}
\label{proCor}
An infinitesimal deformation (\ref{infMulti}) with 
additional series of relations: $t_2^k=0$ for all $k$, is
integrable without any condition on $t_0^k$ and $t_1^k$.
\end{pro}

\begin{proof}
This is an immediate consequence of the fact that all the
Richardson-Nijenhuis products of the two first non-trivial cohomology
classes vanish:
\begin{equation}
\label{ZerProd}
[\![c_0,c_0]\!]=[\![c_0,c_1]\!]=[\![c_1,c_1]\!]=0
\end{equation}
and so do the obstructions.
\end{proof}

The proof that the relations (\ref{condMulti1})-(\ref{condMulti3})
are necessary is just a result of a straightforward computation;
it will be given in Section \ref{proofNec}. 
The proof that these conditions are sufficient is based on
the existence of an important class of deformations corresponding to the
$\Vect(\bbR^n)$-modules of differential operators.

%%%%%%%%%%%%%%%%%%%%%%%%%%%%%%%%%%%%%%%%%%%%%%%%%%%%%%%%%%%%%%%%%%%%%%%%%%%%%%%%%%%%
%%%%%%%%%%%%%%%%%%%%%%%%%%%%%%%%%%%%%%%%%%%%%%%%%%%%%%%%%%%%%%%%%%%%%%%%%%%%%%%%%%%%
\section{Module of differential operators}
%%%%%%%%%%%%%%%%%%%%%%%%%%%%%%%%%%%%%%%%%%%%%%%%%%%%%%%%%%%%%%%%%%%%%%%%%%%%%%%%%%%%
%%%%%%%%%%%%%%%%%%%%%%%%%%%%%%%%%%%%%%%%%%%%%%%%%%%%%%%%%%%%%%%%%%%%%%%%%%%%%%%%%%%%

Consider the space $\cD$ of linear differential operators on $\bbR^n$.
It is isomorphic to $\cS$ as a vector space, but its structure as a
$\Vect(\bbR^n)$-module is quite different. In this section we interpret
$\cD$ as a deformation of the $\Vect(\bbR^n)$-module $\cS$.

%%%%%%%%%%%%%%%%%%%%%%%%%%%%%%%%%%%%%%%%%%%%%%%%%%%%%%%%%%%%%%%%%%%%%%%%%%%%%%%%%%%%
\subsection{Lie derivative of differential operators}\label{ModOp}
%%%%%%%%%%%%%%%%%%%%%%%%%%%%%%%%%%%%%%%%%%%%%%%%%%%%%%%%%%%%%%%%%%%%%%%%%%%%%%%%%%%%

The composition of differential operators is
defined by:
\begin{equation}
A\circ B
=
\sum_{k=0}^{\infty}\frac{1}{k!}\,
\frac{\partial^k A}{\partial \xi_{i_1}\cdots\partial\xi_{i_k}}\,
\frac{\partial^k B}{\partial x^{i_1}\cdots\partial{}x^{i_k}}
\label{compos}
\end{equation}
Of course, since $A$ is a polynomial in $\xi$, there are only
finite number of terms in this sum. 
There is a filtration of the associative algebra $\cD$
\begin{equation}
\label{filtr}
\cD^0
\subset
\cD^1
\subset\cdots\subset
\cD^r
\subset\cdots,
\end{equation}
where $\cD^r$ is the space of $r$-th order differential
operators (isomorphic to $\bigoplus_{i\leq r}\cS_i$
as a vector space).
One has $\cS=\gr\cD$ as well as an associative algebra and as a Lie algebra.
The space $\cS$ is usually called the space of symbols associated to $\cD$.

The space $\cD$ is a $\Vect(\bbR^n)$-module since $\Vect(\bbR^n)$ is
a Lie subalgebra of $\cD$.
Moreover, there is a family of embeddings $\Vect(\bbR^n)\hookrightarrow\cD$
depending on a parameter $\l\in\bbR$ (or $\bbC$) given by
$$
i^\l:X\mapsto
X+\l\,\Div(X)
$$
where $X\in\Vect(\bbR^n)$ and $\Div(X)$ is the divergence
with respect to the standard volume form on $\bbR^n$. This defines
a one-parameter family of $\Vect(\bbR^n)$-module structures on the space $\cD$.
More generally, one can define a two-parameter family of
$\Vect(\bbR^n)$-modules on~$\cD$ by
\begin{equation}
\label{actlm}
\cL_X^{\l,\m}(A)=i^\m(X)\circ A - A\circ i^\l(X)
\end{equation}
These modules are denoted $\cD_{\l,\m}$.

\medskip

\noindent
{\bf Remark}.
From the geometrical viewpoint, the module $\cD_{\l,\m}$ is the space of differential
operators acting on the space of tensor densities (cf. \cite{CMZ,DO,LO2,CM});
the first-order differential operator
$i^\l(X)$ is a Lie derivative of tensor densities of degree $\l$.

\begin{lem}
\label{expForAct}
The explicit formula of the $\Vect(\bbR^n)$-action on $\cD_{\l,\m}$ is
\begin{equation}
\label{actexpl}
\begin{array}{rcl}
\cL_X^{\l,\m} &=& L_X +
(\m-\l)\Div(X)
\\[16pt]
&-&
\displaystyle
\sum_{k=2}^{\infty}
\frac{1}{k!}\left(
\frac{\partial^k X}
{\partial x^{i_1}\!\cdots\!\partial{}x^{i_k}}
\frac{\partial^k}
{\partial \xi_{i_1}\!\cdots\!\partial\xi_{i_k}}
+
k\l
\frac{\partial^{k-1}\Div(X)}
{\partial x^{i_1}\!\cdots\!\partial{}x^{i_{k-1}}}
\frac{\partial^{k-1}}
{\partial \xi_{i_1}\!\cdots\!\partial\xi_{i_{k-1}}}
\right)
\end{array}
\end{equation}
where $L_X$ is as in (\ref{Hamilton}).
\end{lem}
\begin{proof}
This formula readily follows from
(\ref{compos}).
\end{proof}

%%%%%%%%%%%%%%%%%%%%%%%%%%%%%%%%%%%%%%%%%%%%%%%%%%%%%%%%%%%%%%%%%%%%%%%%%%%%%%%%%%%%
\subsection{The Weyl symbols}
%%%%%%%%%%%%%%%%%%%%%%%%%%%%%%%%%%%%%%%%%%%%%%%%%%%%%%%%%%%%%%%%%%%%%%%%%%%%%%%%%%%%

Consider the operator $\Div$ on
$\cD_{\l,\m}$ given by
\begin{equation}
\label{theDiv}
\Div=
\frac{\partial}{\partial x^i}\frac{\partial}{\partial \xi_i}
\end{equation}
that extends the divergence of vector fields to the space of
symmetric contravariant tensor fields. 
Recall that the linear map
\begin{equation}
\label{WeylSym}
\exp\left(
\l\,\Div
\right):\cD\to\cS
\end{equation}
defines the famous
Weyl symbol of a differential operator (see \cite{AW}).
Note that the parameter in this formula is usually interpreted
in terms of the Planck constant, namely $\l=i\hbar/2$.

\begin{lem}
\label{WeylLem}
The action (\ref{actexpl}) becomes after the transformation
(\ref{WeylSym}) as follows:
The action $\widetilde{\cL}^{\l,\m}$ is of the form
\begin{equation}
\label{paramLam}
\widetilde{\cL}_X^{\l,\m} =
L_X+\tau_0c_0(X)+\tau_1c_1(X)+\tau_2c_2(X)+
\sum_{m\geq3}L_m(X)
\end{equation}
with 
\begin{equation}
\label{paramCoef}
\textstyle
\tau_0=\m-\l,
\quad
\tau_1=\l-\half, 
\quad
\tau_2=\l(\l-1).
\end{equation}
where $L_m(X)$ are the terms with the degree shift $m$,
that is, for the operators from $\cS_k$ to $\cS_\ell$ with
$\ell-k=m$.
\end{lem}
\begin{proof}
By definition, $\widetilde{\cL}_X^{\l,\m} =
\exp(-\l\,\Div)\circ\cL_X^{\l,\m}\circ\exp(\l\,\Div)$,
a straightforward computation then yields
(\ref{paramLam}) and (\ref{paramCoef}).
\end{proof}

This new expression of the $\Vect(\bbR^n)$-action on $\cD_{\l,\m}$
allows us to consider this module as a deformation of $\cS$.

%%%%%%%%%%%%%%%%%%%%%%%%%%%%%%%%%%%%%%%%%%%%%%%%%%%%%%%%%%%%%%%%%%%%%%%%%%%%%%%%%%%%
\subsection{Differential operators and formal deformations}\label{DODef}
%%%%%%%%%%%%%%%%%%%%%%%%%%%%%%%%%%%%%%%%%%%%%%%%%%%%%%%%%%%%%%%%%%%%%%%%%%%%%%%%%%%%

The modules $\cD_{\l,\m}$ allow us to prove the existence of a
big class of formal deformations. The idea is to consider
the parameters $\tau_0,\tau_1,\tau_2$ as independent using the fact that
the expressions (\ref{paramCoef}) does not satisfy any non-trivial
homogeneous relation.

\begin{lem}
\label{formaLem}
There exists a (formal) deformation of the form
(\ref{paramLam}) such that the parameters
$\tau_0,\tau_1,\tau_2$ are independent.
\end{lem}

\begin{proof}
Let us use the existence of modules $\cD_{\l,\m}$.
Each term $L_m$ in (\ref{paramLam}) polynomially depends on 
$\tau_0,\tau_1,\tau_2$.
The operator $\widetilde{\cL}_X^{\l,\m}$ defines a 
$\Vect(\bbR^n)$-action and, so, satisfies the homomorphism
condition (\ref{hom}). A term of degree of schift $m$
in (\ref{hom}) is again a polynomial in
$\tau_0,\tau_1,\tau_2$, more precisely, a sum of the terms
\begin{equation}
\label{EqInd}
\tau_0^{m_0}\tau_1^{m_1}\tau_2^{m_2},
\qquad
\hbox{where}
\qquad
m_1+2m_2=m
\end{equation}
with operator coefficients.
But, all the monomials (\ref{EqInd}) with $\tau_0,\tau_1,\tau_2$
given by (\ref{paramCoef}) are, obviously, linearly independent and,
so, the equation (\ref{hom}) has to be satisfied independently
for the operator coefficients
of all monomials (\ref{EqInd}).
These conditions are therefore independent on
$\tau_0,\tau_1,\tau_2$.
\end{proof}

Applying the construction from Section \ref{defGrad} to obtain
a formal deformation with the infinitesimal part
of the form~(\ref{infMulti}),
one then obtains the following intermediate result.

\begin{pro}
\label{intercor}
The following relations:
\begin{equation}
\label{condMulti1'}
t_1^k\,t_2^{k-1} - t_1^{k-2}\,t_2^k =0
\quad,\qquad
k\geq4
\end{equation}
\begin{equation}
\label{condMulti2'}
\quad
\left(t_0^k-t_0^{k-1}\right)
t_1^k
=0
\quad,\qquad
k\geq3
\end{equation}
\begin{equation}
\label{condMulti3'}
\quad
\left(t_0^k-t_0^{k-2}\right)
t_2^k = 0
\quad,\qquad
k\geq4
\end{equation}
are sufficient for 
integrability of the infinitesimal deformation (\ref{infMulti}).
\end{pro}
\begin{proof}
The conditions (\ref{condMulti1'})-(\ref{condMulti3'}) coincide with the
conditions (\ref{relT}) from Proposition \ref{proGen} that are sufficient for
integrability.
\end{proof}

\noindent
{\bf Remark}.
The conditions (\ref{condMulti2'}) and (\ref{condMulti3'}) are 
slicely stronger then  (\ref{condMulti2}) and (\ref{condMulti3})
respectively. So, the ideal generated by these polynomials in
(\ref{condMulti1'})-(\ref{condMulti3'})
is bigger then the one generated by
$R_2(t),R_3(t)$ and $R_3'(t)$.
Therefore, the formal deformation naturally related to the modules
of differential operators turns out to be not the most general one.
In other words, it is not a versal deformation in the sense of
\cite{FF}.

%%%%%%%%%%%%%%%%%%%%%%%%%%%%%%%%%%%%%%%%%%%%%%%%%%%%%%%%%%%%%%%%%%%%%%%%%%%%%%%%%%%%
%%%%%%%%%%%%%%%%%%%%%%%%%%%%%%%%%%%%%%%%%%%%%%%%%%%%%%%%%%%%%%%%%%%%%%%%%%%%%%%%%%%%
\section{Proof of the main theorem}\label{proof}
%%%%%%%%%%%%%%%%%%%%%%%%%%%%%%%%%%%%%%%%%%%%%%%%%%%%%%%%%%%%%%%%%%%%%%%%%%%%%%%%%%%%
%%%%%%%%%%%%%%%%%%%%%%%%%%%%%%%%%%%%%%%%%%%%%%%%%%%%%%%%%%%%%%%%%%%%%%%%%%%%%%%%%%%%

The proof contains two parts. First, we show by a straightforward computation
that the conditions  (\ref{condMulti1}) - (\ref{condMulti3}) are necessary.
Second, we use the existence of the deformation constructed in the preceding
section to prove that these conditions are, indeed, sufficient.

%%%%%%%%%%%%%%%%%%%%%%%%%%%%%%%%%%%%%%%%%%%%%%%%%%%%%%%%%%%%%%%%%%%%%%%%%%%%%%%%%%%%
\subsection{The origin of the integrability conditions}\label{proofNec}
%%%%%%%%%%%%%%%%%%%%%%%%%%%%%%%%%%%%%%%%%%%%%%%%%%%%%%%%%%%%%%%%%%%%%%%%%%%%%%%%%%%%

Let us give here
the details in the case of quadratic relation~(\ref{condMulti1}).

It suffice to look for the solutions of the Maurer-Cartan equation 
which are homogeneous with respect to the partial derivatives in $x$ and $\xi$.
More precisely, one has

\begin{lem}
If there is a solution of the equation (\ref{ordk}), then there exists one
of the form:
\begin{equation}
\label{fromGen}
\begin{array}{rcl}
\vf_m &=&
\displaystyle
\sum_{0\leq{}t\leq{}s\leq3m}
\left(
\a^k_{st}\,
\frac{\partial^{s-t} X}{\partial x^{i_1}\cdots\partial x^{i_{s-t}}}\,
\frac{\partial^{t} }
{\partial x^{i_{s-t+1}}\cdots\partial x^{i_s}}\,
\frac{\partial^{s} }
{\partial \xi_{i_1}\cdots\partial \xi_{i_s}}
\right.\\[16pt]
&&
\displaystyle
\;\;
+\b^k_{st}\,
\frac{\partial^{s-t+1} X}
{\partial x^{i_1}\cdots\partial x^{i_{s-t}}\partial \xi_{i_1}}\,
\frac{\partial^{t} }
{\partial x^{i_{s-t+1}}\cdots\partial x^{i_s}}\,
\frac{\partial^{s-1} }
{\partial \xi_{i_2}\cdots\partial \xi_{i_s}}\\[16pt]
&&
\displaystyle
\;\;
+\left.
\g^k_{st}\,
\frac{\partial^{s-t} X}
{\partial x^{i_1}\cdots\partial x^{i_{s-t-1}}\partial \xi_{i_{s-t}}}\,
\frac{\partial^{t+1} }
{\partial x^{i_{s-t}}\cdots\partial x^{i_s}}\,
\frac{\partial^{s-1} }
{\partial \xi_{i_1}\cdots\partial \widehat{\xi}_{i_{s-t}}
\cdots\partial \xi_{i_s}}
\right)
\end{array}
\end{equation}
\end{lem}

\begin{proof}
The cocycles (\ref{Cocycle0})-(\ref{Cocycle2}) are precisely of this form.
The cup product (\ref{cup}) of two such linear maps is a bilinear map
which is also homogeneous in $x$ and $\xi$. Finally, the coboundary
operator $\d$ preserves the homogeneity in the same way.
\end{proof}

\noindent
{\bf Remark}.
The formula (\ref{fromGen}) is the most general differential operator
on $\cS$ which is invariant with respect to the $\GL(n,\bbR)$-action
on $\bbR^n$.

\medskip

The Maurer-Cartan equation (\ref{ordk}) in the second order
reads:
\begin{equation}
\label{EqMC}
\d\vf_2(t)\left|_{\cS_k}\right.=
-\half\sum_{i,j}t_i^{k-j} t_j^k\,[\![c_i^{k-j},c_j^k]\!]
\end{equation}
Obviously, $[\![c_0,c_0]\!]=[\![c_0,c_1]\!]=0$.
The non-zero cup products are
$$
\begin{array}{rcl}
\half[\![c_1,c_1]\!](X,Y)&=&
-2\,\frac{\partial^2 X}{\partial x^i\partial x^j}\,
\frac{\partial^3 Y}
{\partial x^l\partial x^m\partial \xi_i}\,
\frac{\partial^3 }
{\partial \xi_j\partial \xi_l\partial \xi_m}
-(X\leftrightarrow Y)\\[10pt]
[\![c_0,c_2]\!](X,Y)&=&
3\,\frac{\partial^3 X}{\partial x^i\partial x^j\partial\xi_i}
\frac{\partial^3 Y}{\partial x^l\partial x^m\partial\xi_j}
\frac{\partial^2 }
{\partial \xi_l\partial \xi_m}
-(X\leftrightarrow Y)\\[10pt]
[\![c_1,c_2]\!](X,Y)&=&
-2\,\frac{\partial^2 X}{\partial x^i\partial x^j}
\frac{\partial^4 Y}
{\partial x^l\partial x^m\partial x^p\partial \xi_i}\,
\frac{\partial^4 }
{\partial \xi_j\partial \xi_l\partial \xi_m\partial \xi_p}\\[6pt]
&&+6\,\frac{\partial^3 X}{\partial x^i\partial x^j\partial \xi_l}
\frac{\partial^3 Y}
{\partial x^l\partial x^m\partial x^p}\,
\frac{\partial^4 }
{\partial \xi_i\partial \xi_j\partial \xi_m\partial \xi_p}\\[6pt]
&&
-6\,\frac{\partial^4 X}{\partial x^i\partial x^j\partial x^l \partial \xi_m}
\frac{\partial^3 Y}
{\partial x^m\partial x^p\partial \xi_l}\,
\frac{\partial^3 }
{\partial \xi_i\partial \xi_j\partial \xi_p}\\[6pt]
&&-(X\leftrightarrow Y)\\[10pt]
\half[\![c_2,c_2]\!](X,Y)&=&
-6\,\frac{\partial^3 X}{\partial x^i\partial x^j\partial \xi_l}
\frac{\partial^4 Y}
{\partial x^m\partial x^p\partial x^q\partial \xi_i}\,
\frac{\partial}
{\partial x^l}
\frac{\partial^4 }
{\partial \xi_j\partial \xi_m\partial \xi_p\partial \xi_q}\\[6pt]
&&
9\,\frac{\partial^3 X}{\partial x^i\partial x^j\partial \xi_l}
\frac{\partial^4 Y}
{\partial x^l\partial x^m\partial x^p\partial \xi_q}\,
\frac{\partial}
{\partial x^q}\frac{\partial^4 }
{\partial \xi_i\partial \xi_j\partial \xi_m\partial \xi_p}\\[6pt]
&&
3\,\frac{\partial^3 X}{\partial x^i\partial x^j\partial x^l}
\frac{\partial^4 Y}
{\partial x^m\partial x^p\partial x^q\partial \xi_i}\,
\frac{\partial^5 }
{\partial \xi_j\partial \xi_l\partial \xi_m\partial \xi_p\partial \xi_q}\\[6pt]
&&-3\,\frac{\partial^3 X}{\partial x^i\partial x^j\partial \xi_l}
\frac{\partial^4 Y}
{\partial x^l\partial x^m\partial x^p\partial x^q}\,
\frac{\partial^5 }
{\partial \xi_i\partial \xi_j\partial \xi_m\partial \xi_p\partial \xi_q}\\[6pt]
&&
-6\,\frac{\partial^3 X}{\partial x^i\partial x^j\partial \xi_l}
\frac{\partial^5 Y}
{\partial x^l\partial x^m\partial x^p\partial x^q\partial \xi_i}\,
\frac{\partial^4 }
{\partial \xi_j\partial \xi_m\partial \xi_p\partial \xi_q}\\[6pt]
&&
-(X\leftrightarrow Y)
\end{array}
$$
as well as $[\![c_2,c_0]\!]=[\![c_0,c_2]\!]$
and $[\![c_2,c_1]\!]=[\![c_1,c_2]\!]$.

The second-order term of a formal deformation is:
$\vf_2(t)(X)=\sum{}t_i^k t_j^l\,\vf_{ij}^{kl}(X)$,
where $\vf_{ij}^{kl}(X)$ are differential operators on $\cS$
of the form (\ref{fromGen}) homogeneous with respect to the
partial derivatives in $x$ and in $\xi$ of degree $i+j+2$.
Tedious but direct computation yields:
$$
\begin{array}{rcl}
\vf_2(t)(X) &=&
\a^k_{3}\frac{\partial^3 X}
{\partial x^i\partial x^j\partial x^l}
\frac{\partial^4 }
{\partial \xi_i\partial \xi_j\partial \xi_l}
+\b^k_3\frac{\partial^4 X}
{\partial x^i\partial x^j\partial x^l\partial \xi_i}
\frac{\partial^2 }
{\partial \xi_j\partial \xi_l}
\\[10pt]
&&+
\a^k_{4}\frac{\partial^4 X}
{\partial x^i\partial x^j\partial x^l\partial x^m}
\frac{\partial^4 }
{\partial \xi_i\partial \xi_j\partial \xi_l\partial \xi_m}
+
\g^k_{3}\frac{\partial^4 X}
{\partial x^i\partial x^j\partial x^l\partial \xi_m}
\frac{\partial^4 }
{\partial x^m\partial \xi_i\partial \xi_j\partial \xi_l}\\[10pt]
&&
+\a^k_{5}\frac{\partial^5 X}
{\partial x^i\partial x^j\partial x^l\partial x^m\partial x^p}
\frac{\partial^5 }
{\partial \xi_i\partial \xi_j\partial \xi_l\partial \xi_m\partial \xi_p}
+
\g^k_{4}\frac{\partial^5 X}
{\partial x^i\partial x^j\partial x^l\partial x^m\partial \xi_p}
\frac{\partial^4 }
{\partial x^p\partial \xi_i\partial \xi_j\partial \xi_l\partial \xi_m}
\end{array}
$$
where the coefficients $\a^k_s,\b^k_s,\g^k_s$ are quadratic polynomials
in $t^k_i$ satisfying the following system
\goodbreak
$$
\left\{
\begin{array}{rcl}
3\,\a^k_3 &=& -2\,t_1^{k-1}\,t_1^k\\[4pt]
\b^k_3 &=& 3\,t_0^k\,t_2^k\\[4pt]
4\,\a^k_4+\g^k_3 &=& -2\,t_1^{k-2}\,t_2^k\\[4pt]
2\,\a^k_4 &=& -2t_2^{k-1}\,t_1^k\\[4pt]
\g^k_3 &=& 2\,t_2^{k-1}\,t_1^k\\[4pt]
10\,\a^k_5 &=& -3\,t_2^{k-2}\,t_2^k\\[4pt]
2\,\g^k_4 &=& 3\,t_2^{k-2}\,t_2^k
\end{array}
\right.
$$
This system has a unique solution if and only if the condition
(\ref{condMulti1}) is satisfied. This proves that this condition
is necessary for existence of the second order term $\vf_2(t)$.

The proof in the case of (\ref{condMulti2}) and (\ref{condMulti3}) are
analogous but one has to consider the third-order terms in (\ref{ordk}).

%%%%%%%%%%%%%%%%%%%%%%%%%%%%%%%%%%%%%%%%%%%%%%%%%%%%%%%%%%%%%%%%%%%%%%%%%%%%%%%%%%%%
\subsection{The conditions of integrability are sufficient}\label{proof}
%%%%%%%%%%%%%%%%%%%%%%%%%%%%%%%%%%%%%%%%%%%%%%%%%%%%%%%%%%%%%%%%%%%%%%%%%%%%%%%%%%%%

Let us show that the conditions (\ref{condMulti1})-(\ref{condMulti3}) 
are, indeed, sufficient.

Let us suppose that there is a condition of integrability in order 
$m$, i.e., a relation
$R_m(t)=0$, where $R_m(t)$ is a homogeneous polynomial of degree 
$m$ in $t_0^k,t_1^k,t_2^k$.
One has to prove that the polynomial $R_m(t)$ belongs to the ideal,
$\cR$, generated by
the relations (\ref{condMulti1})-(\ref{condMulti3}).

Proposition~\ref{intercor} insures that $R_m(t)$ belongs to the ideal
generated by the polynomials in (\ref{condMulti1'})-(\ref{condMulti3'}).
Therefore, $R_m(t)$ is split into a sum: 
$R_m(t)=R_{m,1}(t)+R_{m,2}(t)+R_{m,3}(t)$ 
of polynomials divisible by 
(\ref{condMulti1'}), (\ref{condMulti2'}) and (\ref{condMulti3'})
respectively. 

The polynomial $R_{m,1}(t)$
already belongs to $\cR$.

Consider, the second term $R_{m,2}(t)$. 
A direct computation (cf. Section \ref{proofNec}) shows that the only
second-order condition is (\ref{condMulti1}), one then can assume
$m\geq3$. Then, the relation
$[\![c_0,c_1]\!]=0$
implies that each monomial in $R_m(t)$ has to contain some 
parameter~$t_2^\ell$ as a multiple (cf. Proposition~\ref{proCor}).
By assumption, the polynomial $R_{m,2}(t)$ is a multiple of
$\left(t_0^k-t_0^{k-1}\right)t_1^k$ for some $k$.
But, modulo the relation (\ref{condMulti1}), any expression of the form
$\left(t_0^k-t_0^{k-1}\right)
t_1^k\cdots{}t_2^\ell$ is divisible by (\ref{condMulti3}) and so
$R_{m,2}(t)$, indeed, belongs to~$\cR$.

Since the Nijenuis-Richardson product
$[\![c_0,c_2]\!]$ commutes with $c_0$, then
$R_{m,3}(t)$ has to contain the terms of the form
$\left(t_0^k-t_0^{k-2}\right)
t_2^k\cdots{}t_1^\ell$ or 
$\left(t_0^k-t_0^{k-2}\right)
t_2^k\cdots{}t_2^\ell$. But, using the relation
(\ref{condMulti1}) one readily gets that these terms are divisible
by (\ref{condMulti2}) and (\ref{condMulti3}) respectively and, therefore, belong
to $\cR$.

%%%%%%%%%%%%%%%%%%%%%%%%%%%%%%%%%%%%%%%%%%%%%%%%%%%%%%%%%%%%%%%%%%%%%%%%%%%%%%%%%%%%
%%%%%%%%%%%%%%%%%%%%%%%%%%%%%%%%%%%%%%%%%%%%%%%%%%%%%%%%%%%%%%%%%%%%%%%%%%%%%%%%%%%%


\begin{thebibliography}{99}
%%%%%%%%%%%%%%%%%%%%%%%%%%%%%%%%%%%%%%%%%%%%%%%%%%%%%%%%%%%%%%%%%%%%%%%%%%%%%%%%%%%%
%%%%%%%%%%%%%%%%%%%%%%%%%%%%%%%%%%%%%%%%%%%%%%%%%%%%%%%%%%%%%%%%%%%%%%%%%%%%%%%%%%%%

\bibitem{A} 
F. Ammar, {\it Syst\`emes hamiltoniens compl\`etement
integrables et d\'eformations d'alg\`ebres de Lie}.
Publications Math\'ematiques 38 (1994) 427-431.

\bibitem{AW}
G.S.~Argaval, E.~Wolf
{\it Calculus for functions of noncommuting operators and 
general phase space methods in quantum mechanics, I.
Mapping theorems and ordering of functions on noncommuting operators},
Phys. Rev. D, {\bf 2:10} (1970) 2161--2188.

\bibitem{BO} 
S. Bouarroudj, V. Ovsienko,
{\em Three cocycles on ${\rm Diff}(S^1)$
generalizing the Schwarzian Derivative},
Internat. Math. Res. Notices (1998), N.1, 25--39.

\bibitem{CMZ} 
P. Cohen, Yu. Manin, D. Zagier,
{\it Automorphic pseudodifferential operators}, Algebraic
aspects of integrable systems, 17--47, Progr. Nonlinear Differential 
Equations Appl., 26, Birkh\"auser Boston, Boston, MA, 1997.

\bibitem{CM} 
C. Conley, C. Martin
{\it A new family of irreducible representations of the Witt Lie algebra},
to appear in Compositio Math.

\bibitem{DO} 
C. Duval, V. Ovsienko,
{\it Space of second order linear differential operators as a~module over the 
Lie algebra of vector fields}, Advances in Math. {\bf 132}: 2 (1997), 316--333.

\bibitem{FeF} 
 B. L. Feigin, D. B. Fuchs, 
{\it Homology of the Lie algebra of vector
fields on the line}, Func. Anal. Appl., 14 (1980),
201-212.

\bibitem{FF} 
 A. Fialowski, D.B. Fuchs,
{\it Construction of Miniversal Deformations of Lie Algebras},
J. Funct. Anal., {\bf 161} (1999) 76-110.

\bibitem{F} 
D.B. Fuchs, Cohomology of infinite-dimensional Lie
algebras, Consultants Bureau, New York, 1987.

\bibitem{G} 
M. Gersternhaber,
{\it On the deformation of rings and algebras I, III},
Ann. Math. {\bf 79} (1964) 59--103, {\bf 88} (1968) 1--34.

\bibitem{Gar} 
H. Gargoubi, {\it Sur la gŽomŽtrie de l'espace des opŽrateurs diffŽrentiels
linŽaires sur $\bbR$},
Bull. Soc. Roy. Sci. Li\`ge {\bf 69:1} (2000) 21--47.


\bibitem{MP} 
C. Martin, A. Piard, {\it Classification 
of the indecomposable bounded admissible
modules over the Virasoro Lie algebra with weightspaces of 
dimension not exceeding two}, Comm. Math. Phys. {\bf 150}
(1992), no. 3, 465--493.

\bibitem{NR} 
A. Nijenhuis, R.W. Richardson, {\it Deformations of
homomorphisms of Lie algebras}, Bull. AMS {\bf 73} (1967) 175--179.

\bibitem{LO1}
P. Lecomte, V. Ovsienko,
{\it Projectively equivariant symbol calculus}, Lett. Math. Phys.,
{\bf 49:3} (1999) 173--196.

\bibitem{LO2} 
P.B.A. Lecomte, V. Ovsienko,
{\it Cohomology of the vector fields Lie algebra and modules of
differential operators on a smooth
manifold}, Compositio Math. {\bf 119} (2000). 

\bibitem{LN} M. Levy-Nahas, {\it Deformation and contraction 
of Lie algebras.} J. Math. Phys. \textbf{8} (1967) 1211--1222.

\bibitem{MAT}
P.~Mathonet,
{\it Intertwining operators between some spaces of differential operators on a
manifold}, Comm. in Algebra 27 (1999), no. 2, 755--776.

\bibitem{OR1} 
V. Ovsienko, C. Roger,	
{\em Deforming the Lie algebra of vector fields
on~$S^1$ inside the Poisson algebra on $\dot T^*S^1$},
Comm. Math. Phys., 
{\bf 198} (1998) 97--110.


\bibitem{OR} 
V. Ovsienko, C. Roger,
{\it Deforming the Lie algebra of vector fields
on~$S^1$ inside the Lie algebra of pseudodifferential operators on $S^1$},
AMS Transl. Ser.~2, (Adv. Math. Sci.) vol.~194 (1999) 211--227.

\bibitem{R} 
R.W. Richardson,
{\it Deformations of subalgebras
of Lie algebras}, J. Diff. Geom. {\bf 3} (1969) 289--308.



\end{thebibliography}
\end{document}